# The sum of the squares of p positive integers which are consecutive terms of an arithmetic progression: always a non-perfect square when p (a prime) =3 or $p \equiv 5 \pmod{12}$, or $p \equiv 7 \pmod{12}$

*Konstantine Zelator*

*Department of Mathematics*

*317 Bishop Hall*

*SUNY Buffalo State College*

*Buffalo, NY 14222 USA*

*Also: K Zelator*

*P.O Box 4280*

*Pittsburgh, PA 15203 USA*

*E-mail addresses: 1) konstantine-zelator@yahoo.com*

*2) zelatok@buffalostate.edu*

**1. Introduction**

In a paper published by this author in the open access research website, www.academia.edu (see reference [3]), it was established that there exist no three positive integers which are consecutive items of an arithmetic progression; and whose sum of squares is a perfect or integer square. In that paper, we made use of the parametric formulas, which describe all the positive integer solutions of the 4-variable equation,

$$x^2 + y^2 + z^2 = t^2$$

A detailed derivation of those parametric formulas can be found in W. Sierpinski's book,

*Elementary Theory of Numbers* (see reference [2]).

In this work, we offer an alternative proof to this result (that the sum of the squares of three positive integers which are consecutive terms of an arithmetic progression; is always a non-perfect square); a proof that does not make use of the above parametric formulas. This is done in Theorem 1, in Section 2. In Section 3, we state Proposition 1, which states two well known identifies; for the sum of the first k consecutive positive integers; and for the sum of the squares of the first k consecutive positive integers.

In Section 4, we have Propositions 2 and 3. Proposition 2 is a well-known corollary to the Quadratic Reciprocity Law. Namely, that if p and q are distinct odd primes; and at least one of p, q is congruent to 1 modulo 4; then p is a quadratic residue of $p \rightarrow q$ if, and only if, q is a quadratic residue of p.

On the other hand, if $p \equiv q \equiv 3 \ (mod 4)$; then if p is quadratic residue of q; then q is a quadratic non-residue of p. And if p is a non-residue of q; then q is a residue of p.

Proposition 3 is an immediate consequence of Proposition 2. Proposition 3 postulates that if p is a prime such that $p \equiv 5 \ or \ 7 (mod 12)$, then 3 is a quadratic non-residue of p. And if p is a prime such that $p \equiv 1 \ or \ 11 (mod 12)$; then 3 is a quadratic residue of p.

In Section 5, we prove Theorem 2, which postulates that if p is an odd prime and 3 is a quadratic non-residue of p. Then, the sum of the squares of any p positive integers, which are consecutive terms of an arithmetic progression; is a non-perfect square.

Finally, Theorem 3 is an immediate corollary of Theorem 2. It postulates that if p is a prime such that $p \equiv 5 \ or \ 7 (mod 12)$. Then, the sum of the squares of p positive integers, which are consecutive terms of an arithmetic progression; cannot be a perfect square.

## 2. The Case $p = 3$: Theorem 1

**Theorem 1**

*There exist no three positive integers, which are consecutive terms of an arithmetic progression; and whose sum of squares is a perfect or integer square. Equivalently, the sum of the squares of any three positive integers, which are consecutive terms of an arithmetic progression; is a non-perfect square.*

**Proof**

Suppose to the contrary, that there is a triple of three positive integers, which are consecutive terms of an arithmetic progression; and whose sum of squares is a perfect square. We will prove that this assumption leads to a contradiction. It is clear that we may assume that the arithmetic progression is an increasing one; for if it is decreasing, we simply write the three terms in reverse order (also, such a progression could not be constant; for then all three positive integer terms would be equal, which would imply 3 times a perfect square = prefect square, not possible, since $\sqrt{3}$ is an irrational number).

So let $n$ be the first term; among the three consecutive terms; $n + d$ be the second; and $n + 2d$ be the third; where d is the common difference. We have,

$$n^2 + (n + d)^2 + (n + 2d)^2 = t^2 \quad (1)$$

where n, d, and t are positive integers, or equivalently (after expanding the left hand side in equation **(1)**),

$$3n^2 + 6nd + 5d^2 = t^2 \quad \textbf{(2)}$$

$$n, d, t, \in \mathbb{Z}$$

Consider the highest power of 3 diving the positive integers $n, d, t$:

$$n = 3^e * N, \quad d = 3^f * D, \quad t = 3^g * T;$$

where the exponents $e, f, and\ g$ are non-negative integers;

$$e \geq 0, f \geq 0, g \geq 0; e, f, g\ \in \mathbb{Z}$$

And none of the positive integers $N, D, T$; is divisible by 3:

$$3 \nmid N * D * T;\ or\ equivalently\ NDT \not\equiv 0 (mod 3).\ \textbf{(3)}$$

Combining **(1)** with **(2)** yields,

$$3^{2e+1} * N^2 + 2 * 3^{(e+f+1)} N * D + 5 * 3^{2f} * D^2 = 3^{2g} * t^2 \ \textbf{(4)}$$
5

6We will demonstrate that in each case below, a contradiction ensures.

**Case 1:** Assume $e = f$

Then, equation **(4)** takes the form,

$$3^{2e} * [3N^2 + 6ND + 5D^2] = 3^{2g} * T^2 \textbf{ (5)}$$

Since D is not divisible by 3; it is clear that the highest power 3 dividing the left-hand of **(5)** is $3^{2e}$. And the highest power of 3 diving the right-hand side is (by **(3)**) $3^{2g}$. Therefore we must have, $2e = 2g$.

And thus, **(5)** implies that

$$3N^2 + 6ND + 5D^2 = T^2 \textbf{ (6)}$$

Since neither of T and D is divisible by 3; we have,

$$T^2 \equiv D^2 \equiv 1 (mod 3) \textbf{ (7)}$$

Equation **(6)** and congruence **(7)** are contradictory, since they imply that $2 \equiv 1 (mod 3)$, an impossibility. This concludes the first case.

**Case 2:** Assume $e > f$.

If $e > f$; then,

$$2e + 1 > e + f + 1 > 2f \quad (8)$$

Going back to **(4)** and applying **(8)**; we see that

$$3^{2f} * [3^{2(e-f)+1} * N^2 + 2 * 3^{e-f+1} * N * D + 5 * D^2] = 3^{2g} * T^2 \quad (9)$$

Since $3 \nmid D$ and, $(e - f) + 1 \geq 3$ and $e - f + 1 \geq 2$ (by virtue of $e > f$; $e - f \geq 1$, *since* $e, f \in \mathbb{Z}$). Equation **(9)** implies that $2g = 2f$; since $3^{2f}$ and $3^{2g}$ are the highest powers of 3 diving the left and right-hand sides of **(9)** respectively. Cancelling out the common factor $3^{2f} = 3^{2g}$ in **(9)**; leads to the same contradiction modulo 3, as in previous case (Case 1). We omit the details.

**Case 3:** Assume $e < f$.

We distinguish between subcases: the subcase $f = e + 1$; and the subcase $f > e + 1$.

**Subcase 3A:** Assume $f = e + 1$

Then equation **(4)** becomes,

$$3^{2e+1} * [N^2 + 6 * N * D + 15 * D^2] = 3^{2g} * T^2 \quad (10)$$

8Since N is not divisible by 3; equation **(10)** shows that the highest power of 3 dividing the left-hand side is $3^{2e+1}$; while $3^{2g}$ is the highest power of 3 diving the right-hand side; which renders **(10)** contradictory since obviously $2e + 1$ cannot equal $2g$.

**Subcase 3b:** Assume $f > e + 1$.

This assumption implies,

$$2e + 1 < e + f + 1 < 2f \quad \textbf{(11)}$$

Going back to equation **(4)** and applying **(11)** implies,

$$3^{2e+1} * \left[N^2 + 2 * 3^{f-e} * N * D + 5 * 3^{2(f-e)} * D^2\right] = 3^{2g} * T^2 \quad \textbf{(12)}$$

Since the exponents $(e - f)$ and $2(f - e)$ are positive integers; and since N is not divisible by 3. We obtain the same contradiction as in the previous subcase: $3^{2e+1}$ exactly diving the left-hand side of **(12)**; while $3^{2g}$ exactly dividing the right-hand side. The proof is complete. ∎

## 3. Two Well-Known Identities

The following two identities are very well known; and may be stated without proof.

**Proposition 1**

Let k be a positive integer. Then,

$$(i) \sum_{i=1}^{k} = 1 + 2 = \cdots + k = \frac{k(k+1)}{2}$$

$$(ii) \sum_{i=1}^{k} i^2 = 1^2 + 2^2 + \cdots + k^2 = \frac{k(k+1)(2k+1)}{6}$$

## 4. Some basic information on quadratic reciprocity and quadratic residues. Also, Propositions 2 and 3.

The following information is very well-known, and can be readily found in number theory books; for example, see references [1] and [4].

Let p be an odd prime; and $a$ an integer not divisible by p. Then $a$ is called a *quadratic residue of p* if the one-variable congruence, $x^2 \equiv a(mod p)$; has a solution. It can be shown, that if this congruence has a solution it has *exactly two* incongruent solutions modulo p.

As it turns out; among the $(p-1)$ integers, $1, 2, \ldots, p-1$; exactly half are quadratic residues of p; the other $\left(\frac{p-1}{2}\right)$ among them are not quadratic residues of p; they are called *quadratic non-residues* of p. Furthermore, it can be shown that if an integer a (not divisible by p) is a quadratic residue of p; then a must be a congruent of modulo p; to one of the squares: $1^2, 2^2, \ldots, \left(\frac{p-1}{2}\right)^2$.

**The Legendre Symbol:** *If a is not divisible by p; then the notation $\left(\frac{a}{p}\right) = 1$; means that a is a quadratic residue of p; while $\left(\frac{a}{p}\right) = -1$ means that a is a quadratic non-residue of p.*

**Quadratic Reciprocity Law:** *Let p and q be distinct odd primes. Then,*

$$\left(\frac{p}{q}\right)\left(\frac{q}{p}\right) = (-1)^{\left(\frac{p-1}{2}\right)\left(\frac{q-1}{2}\right)}.$$

An immediate Corollary to the Quadratic Reciprocity Law is the following:

*If p and q are distinct odd primes then,*

$$\left(\frac{p}{q}\right)\left(\frac{q}{p}\right) = 1 \; if \; p \equiv 1 (mod 4) \, or \, q \equiv 1 (mod 4); \; or \; -1 \; if \; p \equiv q \equiv 3 (mod 4)$$

We rewrite the previous Corollary, by using words, instead of the Legendre symbols.

**Proposition 2**

*Suppose that p and q are distinct odd primes. Then,*

*(a) If at least one p, q; is congruent to 1 modulo 4:*

*Either p and q are quadratic residues of each other or, they are quadratic non-residues of each other.*

*(b) If $p \equiv q \equiv 3 \pmod 4$. Then p is a quadratic residue of q; and q is a quadratic non-residue of p (or vice versa).*

Now consider a prime p modulo 12. Clearly any positive integer greater than 3 and congruent to $2, 3, 4, 6, 8, 9, \text{ or } 10 \pmod{12}$; must be a composite (i.e. nonprime) number. Thus if p is a prime then: $p \equiv 1, 5, 7, \text{ or } 11 \pmod{12}$.

If $p \equiv 1 \pmod{12}$; then $p \equiv 1 \pmod 3$ and $p \equiv 1 \pmod 4$.

But then p is a quadratic residue of 3. And thus, by Proposition 2(a), it follows that 3 is a quadratic residue of p. If $p \equiv 5 \pmod{12}$; then $p \equiv 1 \pmod 4$ and $p \equiv 2 \pmod 3$. Using the fact that 2 is a quadratic non-residue of 3, and Proposition 2(a); we see that in this case, 3 is a quadratic non-residue of p. We use Proposition 2(b) and a similar argument for the cases $p \equiv 7 \pmod{12}$; and $p \equiv 11 \pmod{12}$. We have the following.

**Proposition 3**

Let p be a prime, $p \geq 5$. Then

$(i)$ If $p \equiv 5 \pmod{12}, 3 \text{ is a quadratic non} - \text{residue of } p.$

$(ii)$ If $p \equiv 7 \pmod{12}, 3 \text{ is a quadratic non} - \text{residue of } p.$

$(iii)$ If $p \equiv 1 \pmod{12}, \text{ or if } p \equiv 11 \pmod{12}.$ Then 3 is a quadratic residue of p.

## 5. Theorems 2 and 3

**Theorem 2**

*Let p be an odd prime with $p \geq 5$. And suppose that 3 is a quadratic non-residue of p. Then, there exist no p positive integers, which are consecutive terms of an arithmetic progression and whose sum of squares is a perfect square. Equivalently, the sum of the squares of any p positive integers, which are consecutive terms of an arithmetic progression; is a non perfect square.*

**Proof**

We argue by contradiction. Suppose to the contrary that,

$$n^2 + (n+d)^2 + \cdots [n + (p-1)d]^2 = t^2 \text{ for positive integers n, d, and t.} \quad \textbf{(13)}$$

By expanding the left-hand side of equation **(13)** we obtain,

$$p*n^2 + 2[1 + 2 + \cdots + (p-1)]nd + [1^2 + 2^2 + \cdots + (p-1)^2]d^2 = t^2 \quad \textbf{(14)}$$

Next, we make use in **(14)** of the two identities in Proposition 1 to get,

$$p*n^2 + p*n*d(p-1) + \frac{p(p-1)(2p-1)}{6} = t\text{\textasciicircum}2; \text{ or } 6pn^2 + 6p(p-1)nd +$$
$$p(p-1)(2p-1)d^2 = 6t^2. \quad \textbf{(15)}$$

Let $p^e, p^f, p^g$; be the highest powers of the prime dividing the integers n, d, and t respectively:

$$n = p^e * N, d = p^f * D, t = p^g * T$$

where the exponents e, f, g, are non-negative integers; $e \geq 0, f \geq 0, g \geq 0$; $e, f, g \in \mathbb{Z}$. And none of the positive integers N, D, T is **(16)** divisible by p: $p$ does not divide $N * D * T$; or equivalently $NDT \not\equiv 0 \pmod{p}$.

Combining **(16)** and **(17)** yields,

$$6 * p^{2e+1} * N^2 + 6 * p^{e+f+1}(p-1)N * D + p^{2f+1} * (p-1)(2p-1) * D^2 = 6 * p^{2g} * T^2 \quad \textbf{(17)}$$

First, observe that if $e > f$ or $e < f$; we obtain and obvious contradiction. Indeed, if $e > f$, then equation **(17)** implies (since then $2e + 1 > e + f + 1 > 2f + 1$),

$$p^{2f+1} * \left[ 6 * p^{2(e-f)} * N^2 + 6 * p^{(e-f)}(p-1)N * D + (p-1)(2p-1) * D^2 \right] = 6p^{2g} * T^2 \quad \textbf{(18)}$$

Equation **(18)** is contradictory. Indeed, since D is not divisible by p (by **(16)**) and obviously neither of the $p - 1, 2p - 1$ is divisible by p; it follows (since p is a prime) that the product $(p - 1)(2p - 1) * D^2$ is not divisible by p. Therefore, by virtue of $e - f \geq 1$; it follows that the highest power of p dividing the left-hand side of **(18)**, is $p^{2f+1}$; which obviously cannot equal $p^{2g}$, the highest power of p dividing the right-hand side of (18) (since by (16), T is not divisible by p). A similar argument can be given in the case $e < f$; we omit the details.

**The case $e = f$**

This is the only remaining case. Accordingly, if $e = f$, then equation **(17)** takes the form,

$$p^{2e+1} * [6N^2 + 6(p-1)ND + (p-1)(2p-1)d^2] = 6 * p^{2g}T^2 \quad \textbf{(19)}$$

Obviously, $2e + 1 \neq 2g$. It follows then from equation **(19)** that $2g > 2e + 1$. Indeed, the possibility $2g < 2e + 1$ is ruled out; for, it would imply that the prime p divides T; contrary to the conditions in **(16)**. Therefore, we must have $2g > 2e + 1$; which by **(19)** implies that,

$$6N^2 - 6ND + 6pND + (p-1)(2p-1)D\wedge 2 \equiv 0 (mod p). \quad \textbf{(20)}$$

Also $(p-1)(2p-1)D^2 = 2p^2D^2 - pD^2 - 2pD^2 + D^2$; which implies

$$(p-1)(2p-1)D^2 \equiv D\wedge 2 (mod p) \quad \textbf{(21)}$$

From **(20)** and **(21)** we further obtain that

$$6N^2 - 6ND + D^2 \equiv 0 (mod\, p); \text{ or equivalently, } (3N - D)^2 \equiv 3N^2 (mod\, p). \quad \textbf{(22)}$$

Recall that neither of N, D is divisible by the prime $p \geq 5$. Which in turn implies by **(22)**; that neither of the integers $3N - D$ and $N$; is divisible by p. Since $N \not\equiv 0 \ (mod\, p)$; there exist a unique modulo p, integer K such that,

$$N * K \equiv 1 (mod\, p). \quad \textbf{(23)}$$

(The integer K is often referred to as the multiplicative inverse of N modulo p.) By multiplying both sides of (22) by $K^2$, we obtain $K^2 * (3N - D)^2 \equiv 3 * N^2 K^2 (mod\, p)$; and since $N^2 K^2 \equiv 1 (mod\, p), by\ (23)$.

We arrive at,

$[K(3N - D)]^2 \equiv 3 (mod\, p)$, which implies that 3 is a quadratic non-residue of p; contrary to the hypothesis of the theorem. The proof is complete. ∎

The following theorem, Theorem 3, is an immediate consequence of Theorem 2 and Proposition 3.

**Theorem 3**

*Let p be an odd prime, $p \geq 5$, and $p \equiv 5$ or $7 \pmod{12}$. Then, there exist no p positive integers, which are consecutive terms of an arithmetic progression; and whose sum of squares is a perfect square. Equivalently, the sum of the squares of any p positive integers, which are consecutive terms of an arithmetic progression; is a non-perfect square.*